\documentclass[12pt ,twoside ,a4paper]{article}
\usepackage{enumitem}
\setlist{noitemsep, topsep=0pt, parsep=0pt, partopsep=0pt, leftmargin=0.5cm}
\normalfont
\usepackage{multirow}
\usepackage{hyperref}
\usepackage{natbib}
\usepackage{amsthm}
\newtheorem{theorem}{Theorem}[section]
\newtheorem{lemma}[theorem]{Lemma} 
\newtheorem{remark}{Remark}[section]
\newtheorem{definition}[theorem]{Definition}

\usepackage[margin=1in]{geometry}
\usepackage{amsmath}

\usepackage[utf8]{inputenc}
\usepackage{geometry}
\usepackage{amssymb}
\usepackage{breqn}

\title{\textbf{On asymptotic behavior of U-statistics for associated random variables}}
\vspace{-2.0pt}
\author{${\text{ Mansi Garg}\footnote{Corresponding author.}}$ and $\text{Isha Dewan}$ \vspace{-0.1in}\\
 ${}$Indian Statistical Institute \vspace{-0.1in}\\New Delhi-110016 (India)\vspace{-0.1in}\\ mansibirla@gmail.com and ishadewan@gmail.com $\vspace{-0.10in}$ }
 
\begin{document}
\fontsize{15}{20}
\date{}
\maketitle
\abstract{Let $\lbrace{X_n, n \ge 1}\rbrace$ be a sequence of stationary associated random variables. For such a sequence we discuss the limiting behavior of U-statistics based on kernels which are of bounded Hardy-Krause variation. }\\
\\
\textbf{Keywords}: {\it Associated random variables; Central limit theorem; Gini's mean difference;  Hardy-Krause variation;  U-statistics.}

\numberwithin{equation}{section}
\section{Introduction}

  In this paper we discuss the asymptotic behavior of U-statistics which are based on kernels of bounded Hardy-Krause variation. Apropos our discussion, we give the following definition.

\begin{definition}(Esary et al. (1967))  A finite collection of random variables $\{X_j,  1 \le j \le n \}$  is said to be associated, if for any choice of component-wise nondecreasing functions $h$, $g$ $:$ ${\mathbb{R}}^n \rightarrow {\mathbb{R}}$, we have,
\begin{equation*}
Cov(h(X_1,\ldots, X_n), g(X_1,\ldots, X_n)) \ge 0
\end{equation*}
whenever it exists. An infinite collection of random variables $\{X_j,  j \ge 1 \}$ is associated if every finite sub-collection is associated. 
\end{definition}

Associated random variables have been widely used in reliability studies, statistical mechanics, and percolation theory. A set consisting of independent random variables is associated (cf. Esary et al. $(1967)$). Monotonic functions of associated random variables are associated. For details on probabilistic results and examples relating to associated sequences, see Bulinski and Shashkin $(2007, 2009)$, Prakasa Rao $(2012)$, and Oliveira $(2012)$.

Given stationary associated observations $\lbrace{X_j, 1 \le j \le n}\rbrace$, the U-statistic $U_n$ of degree $k$ $( 1\le k \le n)$ based on a symmetric kernel $\rho: \mathbb{R}^k \to \mathbb{R}$ is given by,
\hypertarget{definition 1.1}{\begin{equation}
U_n = {{n}\choose{k}}^{-1}\sum_{(n, k)} {\rho}(X_{j_1}, X_{j_2}, ..., X_{j_k}),
\end{equation}}
where $(n, k)$ indicates all subsets $1\le {j_1}< {j_2}< ...< {j_k}\le n$ of ${\left\lbrace1, 2, \dots, n\right\rbrace}$. 

Let $F$ be the distribution function of $X_1$ and $U_n$ be the U-statistic based on the symmetric kernel ${\rho}(x_1, x_2)$. Let $\theta = \int\limits_{\mathbb{R}^2}{\rho}(x_1, x_2)~{dF({x_1})dF({x_2})}$.  Define, 
\begin{align*}
&  {\rho}_{1}({x_1}) = \int\limits_{\mathbb{R}}{\rho}(x_1, x_2)~dF({x_{2}}), \: \:  h^{(1)}({x_1}) =  {\rho}_1({x_1}) - \theta,\\
&  \text{and }   \: h^{(2)}({x_1}, {x_2}) =  {\rho}({x_1}, {x_2}) - {\rho}_1({x_1}) - {\rho}_1({x_2}) + \theta.
  \end{align*}
   Then, the Hoeffding-decomposition (H-decomposition) for $U_n$ is (see, Lee $(1990)$) $U_n = \theta + 2H_n^{(1)} + H_n^{(2)}$, where $H_{n}^{(j)}$ is the U-statistic of degree $j$ based on the kernel $h^{(j)}$, $j = 1, 2$. When the observations are $i.i.d$, $E(U_n) = \theta$.

Next, we discuss the concepts of Hardy-Krause variation and Vitali variation. Discussions and applications of these concepts can be found in   Clarkson and Adams $(1933)$, Adams and Clarkson $(1934)$, Owen $(2004)$ and Beare $(2009)$. The following is from Beare $(2009)$.

\begin{definition}  The Vitali variation of a function $f: [a, b] \rightarrow \mathbb{R}$, where $[a, b] = \lbrace{\textbf{x} \in \mathbb{R}^k : a \le \textbf{x} \le b}\rbrace$, $a, b \in \mathbb{R}^k$, $k \in \mathbb{N}$  is defined as
$||f||_V = sup \sum_{R \in \mathbb{A}} |\Delta_R f|$. The supremum is taken over all finite collections of k-dimensional rectangles $\mathbb{A} = \lbrace{R_i : 1 \le i \le m}\rbrace$ such that $\cup_{i=1}^{m}R_i = [a, b]$, and the interiors of any two rectangles in $\mathbb{A}$ are disjoint. Here, if $R = [c, d]$, a k-dimensional rectangle contained in $[a, b]$, then, $\Delta_R f= \sum_{I \subseteq \lbrace {1, 2,..., k}\rbrace} (-1)^{|I|}f(x_I)$, where, $x_I$ is the vector in $\mathbb{R}^k$ whose $i^{th}$ element is given by $c_i$ if $i \in I$, or by $d_i$ if $i \not \in I$, $f_\emptyset$ $=$ $f(b)$. For instance, if $k=2$ and $R = [c_1, d_1]\times [c_2, d_2]$ then, $\Delta_Rf = f(d_1, d_2) - f(c_1, d_2)  - f(d_1, c_2)  + f(c_1, c_2)$.
\end{definition}

\begin{definition} The Hardy-Krause variation of a function $f: [a, b] \rightarrow \mathbb{R}$, $[a, b] = \lbrace{\textbf{x} \in \mathbb{R}^k : a \le \textbf{x} \le b}\rbrace$, $a, b \in \mathbb{R}^k$, $k \in \mathbb{N}$, is given by,
$||f||_{HK} = \sum_{\emptyset \neq I \subseteq \lbrace{1, ..., k}\rbrace} ||f_I||_V$. Here, given a non-empty set $I \subseteq \lbrace{1, 2, ..., k}\rbrace$, $f_I$ denotes the real valued function on $\prod_{i \in I}[a_i, b_i]$ obtained by setting the $i^{th}$ argument of $f$ equal to $b_i$ whenever $i \not \in I$.
\end{definition}
When $k=1$, the Hardy-Krause variation is equivalent to Vitali variation and hence the standard definition of total variation. 

 If  $f: [a, b] \rightarrow \mathbb{R}$, $[a, b] = \lbrace{\textbf{x} \in \mathbb{R}^k : a \le \textbf{x} \le b}\rbrace$, $a, b \in \mathbb{R}^k$, $k \in \mathbb{N}$ is of bounded Hardy-Krause variation, then for any $x \in (a, b]$ there exists a value, denoted by $f^{-}(x)$ such that $f(x_m) \to f^{-}(x)$ for any sequence of points $\lbrace{x_m}\rbrace$ $ {\in}$ $[a, x)$ that converges to $x$. Set $f^{-}(x) = f(x)$ for $x \not {\in} (a, b]$. $f^{-}$ is referred as the left-hand limit of $f$. If $f_I = f_I^-$ for all non-empty $I \subseteq \lbrace{1, ..., n}\rbrace$, then we say $f$ is left-continuous. Similarly, the right-hand limit $f^{+}$ and the right-continuity of $f$ can be defined.

 The paper is organized as follows. Section 2 includes results and definitions that will be required to prove our main results in section 3. In section 3 of the paper, we obtain a central limit theorem for U-statistics based on functions of bounded Hardy-Krause variation for stationary associated random variables.  In section 4, we apply our results to obtain the asymptotic distribution of Gini's mean difference. We give simulation results in section 5 to investigate the asymptotic normality of the statistic under the dependent setup.

\section{Preliminaries} 
\setlength\parindent{24pt}
 In this section, we give results and definitions which will be needed to  prove our main results given in section 3.

{\begin{lemma}
(Newman (1980)) Let $\left\lbrace X_n,  n \ge 1\right\rbrace$ be a strictly stationary sequence of associated random variables. Let $\sigma^{2}$ $=$ $Var( X_1)$ $+$ $2\sum_{j =2}^{\infty}$ $Cov( {X_1}, {X_j})$ with $0< \sigma^{2}< \infty$. Then,
\begin{equation}
\frac{1}{\sqrt{n}\sigma}\sum_{j=1}^{n}(X_j - E(X_j)) \xrightarrow {\mathcal{ L}} N(0, 1) \: \text{as}\: {n \to \infty}.
\end{equation}
\end{lemma}}

{\begin{definition} (Newman (1984))  Let $f$ and $f_1$ be two complex-valued functions on $\mathbb{R}^{n}$, $n \in \mathbb{N}$. We say $ f \ll f_1$ if $f_1 - Re(e^{i\alpha}f)$ is coordinate-wise nondecreasing for every $\alpha \in \mathbb{R}$. If $f$ and $f_1$ are two real-valued functions on $\mathbb{R}^{n}$, then $ f \ll f_1$ iff $f_1 + f$ and $f_1 - f$ are both coordinate-wise nondecreasing. If $f \ll f_1$, then $f_1$ will be coordinate-wise nondecreasing. 
\end{definition}}

{\begin{lemma}
(Newman (1984) Let $\left\lbrace X_n, n \ge 1\right\rbrace$  be a stationary sequence of associated random variables. For each j, let $Y_j = f(X_j)$ and $\tilde{Y_j} = \tilde{f}(X_j)$. Suppose that $f \ll \tilde{f}$. Define $\sigma^2 = Var (Y_1) + 2\sum_{j=2}^{\infty}Cov(Y_1, Y_j)$. Let $\sigma^2 > 0$ and $0 < \sum_{j=1}^{\infty}Cov(\tilde{Y_1},\tilde{Y_j}) < \infty$.
Then,
\begin{equation}
\frac{1}{\sqrt{n}\sigma}\sum_{j=1}^{n}(Y_j - E(Y_j)) \xrightarrow {\mathcal{ L}} N(0, 1) \: \text{as}\: {n \to \infty},
\end{equation}
\end{lemma}}

{\begin{lemma}
(Lebowitz (1972)) Let the random variables $\lbrace{X_j, 1 \le j \le n}\rbrace$ be associated. Define, for $A$ and $B$ subsets of $\lbrace{1, 2, ..., n}\rbrace$ and real $x_j `s$,
\begin{equation*}
H_{A, B} (x_j, j \in A \cup B) = P[ X_j > x_j; j \in A \cup B] - P[X_k > x_k, k \in A]P[X_l > x_l, l \in B].
\end{equation*}
Then,
\begin{equation*}
0 \le H_{A, B} \le \sum_{i \in A}\sum_{j \in B} H_{\lbrace{i}\rbrace, \lbrace{j}\rbrace}.
\end{equation*}

\end{lemma}}

{\begin{lemma}
(Demichev (2014)) Suppose $X$ and $Y$ are associated random variables with bounded continuous densities and $(X, Y) \in L_2$. Then,  for any $0 < \delta < 1/2$ there exists $C = C(\delta)$ such that,
\begin{equation}
\underset{x,y}{sup} |P( X \le x, Y \le y) - P( X \le x)P(Y \le y)| \le C[Cov(X, Y)]^{\delta}.
\end{equation}
\end{lemma}}

\begin{lemma} (Beare (2009))
Let  $Z$ be a random vector taking values in a bounded $u+v$-dimensional rectangle $R = [a, b]$ $\subset$ $\mathbb{R}^{u+v}$, $u, v \in \mathbb{N}$. $R$ is chosen such that each $Z_i$ is equal to $a_i$ with probability zero. For a non-empty set $K \subseteq \lbrace{1, ..., u+v}\rbrace$, let $R_K$ $=$ $\prod_{k \in K}[a_k, b_k]$, and let $F_K$ denote the joint distribution of those $Z_k$ for which $k \in K$. Let $F_\emptyset = 1$. $\textbf{X}=(Z_1, ..., Z_u)$ and $\textbf{Y}=(Z_{u+1}, ..., Z_{u+v})$ and the two real functions $f$ and $g$ are defined on $R_{\lbrace{1, ..., u}\rbrace}$ and $R_{\lbrace{u+1, ..., u+v}\rbrace}$. Suppose $f$ and $g$ are of bounded Hardy-Krause variation ($||f||_{HK}, ||g||_{HK} $ $<$ $\infty$) and left-continuous.  If $\gamma$ $<$ $\infty$ is such that, $|| F_{I \cup J} - F_{I}F_{J}||_{\infty} \le \gamma$, for all non-empty sets $I \subseteq \lbrace {1, 2,..., u}\rbrace$ and $J \subseteq \lbrace {u+1, u+2,..., u+v}\rbrace$, then we have,
\begin{equation*}
|Cov(f(\textbf{X}), g(\textbf{Y}))| \le \gamma ||f||_{HK}||g||_{HK}.
\end{equation*}
\end{lemma}

{\begin{lemma} (Garg and Dewan (2015)) Let $\left\lbrace X_n, n \ge 1\right\rbrace$  be a stationary sequence of associated random variables. For each j, let $Y_j = f(X_j)$ and $\tilde{Y_j} = \tilde{f}(X_j)$ such that $f \ll \tilde{f}$. Let $E(Y_1) = \mu$ and $0 < E(Y^2_1) < \infty$. Suppose $\left\lbrace \ell_n, n \ge 1 \right\rbrace$  is a sequence of positive integers with $1\le \ell_n\le n$ and  $\ell_n$ $=$ $o(n)$ as $ n \to \infty$. Set $S_j(k) = \sum_{i= j +1}^{j+ k}Y_i$, $\bar{Y_n} =\frac{1}{n} \sum_{j =1}^{n} Y_i$.  Define,  $(write\hspace{2 mm} \ell = \ell_n)$,
\begin{equation}
B_n = \frac{1}{n-\ell+1} \Bigg(\sum_{j=0}^{n-\ell} \frac{| S_j(\ell) - \ell\bar{Y_n}|}{\sqrt{\ell}}\Bigg).
\end{equation}
Assume, $\sum_{j=1}^{\infty}Cov(\tilde{Y_1},\tilde{Y_j}) < \infty$. Then, 
\begin{equation*}
B_n \to \sigma_f\sqrt{\frac{2}{\pi}}  \hspace{3 mm} \text{in}  \hspace{2 mm} L_2 \hspace{2 mm} \text{as} \hspace{2 mm}n \to \infty,
\end{equation*}
where, $\sigma_f^2 = Var(Y_1) + 2\sum_{j = 2}^{\infty} Cov (Y_1, Y_j)$.
\end{lemma}}

{\begin{lemma} (Birkel (1988a))  Let $\left\lbrace X_n, n \ge 1\right\rbrace$  be a stationary sequence of associated random variables with $E(X_j) = 0$ and $|X_j| \: \le C_1 < \infty$ for $j \ge 1$.
Assume that $u(n) = 2\sum_{j=n+1}^{\infty} Cov(X_1, X_j) = O( n^{-(r-2)/2})$. Then, there is a constant $B > 0$ not depending on $n$ such that for all $n \ge 1$,
\begin{equation}
\underset{m \ge 0}{sup} E|S_{n+m} - S_m|^r \le Bn^{r/2},
\end{equation}
where, $S_n = \sum_{j=1}^n X_j$.
\end{lemma}}

\section{Limiting behavior of U-statistics based on kernels of bounded Hardy-Krause variation.}

The main result of this section is Theorem {$3.3$}. It gives the central limit theorem for U-statistics based on a kernel of degree 2 which is of bounded Hardy-Krause variation given a sequence of stationary associated random variables. The extension of this theorem to U-statistics with kernels of a general finite degree k $\ge$ $3$ is also discussed. We also discuss a strong law of large numbers for U-statistics based on such kernels using the results in Christofides $(2004)$.

{\subsection{Central limit theorem}
\begin{lemma} 
Let $\left\lbrace X_n, n \ge 1 \right\rbrace$ be a sequence of stationary associated random variables with $|X_n| < C_1 < \infty$, for all $n \ge 1$. Assume that the density function of $X_1$, denoted by $f$, is bounded.  If $h^{(2)}(x, y)$ is a degenerate symmetric kernel of degree $2$ which is of  bounded Hardy-Krause variation and left continuous, then, under the condition $\sum_{j=1}^\infty Cov(X_1, X_j)^{\gamma}$ $<$ $\infty$, for some $0 < \gamma< 1/6$,
\begin{equation}
\sum_{1 \le i < j \le n} \sum_{1 \le k < l \le n} |E( h^{(2)}(X_i, X_j) h^{(2)}(X_k, X_l))| = o(n^3).
\end{equation}
\begin{proof} Let $C$ be a generic positive constant in the sequel. Let $\lbrace{X_i', i \ge 1}\rbrace$ be a sequence of random variables independent of $\lbrace{X_i, i \ge 1}\rbrace$ such that $\lbrace{X_i', i \ge 1}\rbrace$ are  $i.i.d$ with $f$ as the marginal density function of $X_1'$. Observe that from definition $h^{(2)}(x, y)$ is a degenerate kernel, i.e. $\int_{\mathbb{R}} h^{(2)}(x, y)dF(x) = 0$ for all $y \in \mathbb{R}$. Hence,
\begin{equation}
E( h^{(2)}(X_i', X_j) h^{(2)}(X_k, X_l)) = 0.
\end{equation}
 Let $K = \lbrace{i, j, k, l}\rbrace$, and for all $\emptyset \neq A \subseteq K$, define $I= A \cap \lbrace{i}\rbrace$ and $J= A \cap \lbrace{j, k, l}\rbrace$ ($A = I \cup J$). Let the joint distribution function of $\lbrace{X_a, a \in S}\rbrace$ be denoted by $F_S$, for any $S \subseteq A$ ($F_{\emptyset} = 1$). Define $\tilde{F}_A=  F_{A} - F_ IF_J$. For any  $\textbf{x} = (x_i, x_j, x_k, x_l) \in [-C_1, C_1]^4$, define $ z(\textbf{x}) = h^{(2)}(x_i, x_j)h^{(2)}(x_k, x_l)$, and $z_A(\textbf{x}) = z(\textbf{x}_A)$ where, $\textbf{x}_A$ is obtained by setting $x_t$ in $\textbf{x}$ equal to $C_1$ whenever $t \notin A$.  $z_A$ is then a real valued function on $[-C_1, C_1]^{|A|}$. Observe that $z_A$ is also a function of bounded Hardy-Krause variation and left-continuous. Let $\mu_{z_A}$ be the signed measure generated by $z_A$. (A discussion on the construction of measures from functions of bounded Hardy-Krause variation can be found in Beare $(2009)$ and the references therein.) Replicating the arguments of Theorem $3.1$ of Beare $(2009)$, for all $1 \le i < j \le n$ and $1 \le k < l \le n$, we get,
\begin{align}
& \Big|E(z(X_i, X_j, X_k, X_l))\Big| =   \Big| E(z(X_i, X_j, X_k, X_l)) - E(z(X'_i, X_j, X_k, X_l))\Big|  \nonumber \\
 & = \Big| \sum_{I = \lbrace{i}\rbrace, \emptyset \neq J \subseteq  \lbrace{j, k, l}\rbrace} (-1)^{|I \cup J|}\int_{[-C_1, C_1]^{|I \cup J|}} \tilde{F}_{I \cup J} d\mu_{z_{I \cup J}}\Big| \\
 & \le \sum_{\emptyset \neq J \subseteq  \lbrace{j, k, l}\rbrace}  ||  \tilde{F}_{\lbrace{i}\rbrace \cup J}||_{\infty} (|\mu_{z_{\lbrace{i}\rbrace \cup J}}|[-C_1, C_1]^{|\lbrace{i}\rbrace \cup J|}).
\end{align}
 where, $|\mu_{z_{\lbrace{i}\rbrace \cup J}}|$ denotes the total variation of the measure $\mu_{z_{\lbrace{i}\rbrace \cup J}}$. Following the ideas illustrated in Theorem $4.1$ of  Beare $(2009)$, $(|\mu_{z_A}|[-C_1, C_1]^{|A|} ) = ||z_A||_V$, where, $||z_A||_V$ denotes the Vitali variation of $z_A$. Since a function that is of bounded Hardy-Krause variation is of bounded Vitali variation, 
\begin{align}
 |E( h^{(2)}(X_i, X_j) h^{(2)}(X_k, X_l))|  & \le C\sum_{\emptyset \neq J \subseteq  \lbrace{j, k, l}\rbrace}  ||  \tilde{F}_{\lbrace{i}\rbrace \cup J}||_{\infty} \nonumber \\
& \le C(Cov(X_i, X_j)^{\delta} + Cov(X_i, X_k)^{\delta}  + Cov(X_i, X_l)^{\delta} ), 
\end{align}
for  $0 < \delta = 3\gamma  < 1/2$. The last inequality in $(3.5)$ follows from Lemmas $2.4$ and $2.5$.
Similarly, 
\begin{align}
& |E( h^{(2)}(X_i, X_j) h^{(2)}(X_k, X_l))|  \nonumber \\
& \le C(Cov(X_j, X_i)^{\delta} + Cov(X_j, X_k)^{\delta}  + Cov(X_j, X_l)^{\delta} ), \\
& |E( h^{(2)}(X_i, X_j) h^{(2)}(X_k, X_l))| \nonumber \\
& \le C(Cov(X_k, X_j)^{\delta} + Cov(X_k, X_i)^{\delta}  + Cov(X_k, X_l)^{\delta} ), \text{for  $0 < \delta = 3\gamma < 1/2$}.
\end{align}
Combining $(3.5)$, $(3.6)$ and $(3.7)$,
\begin{equation}
 |E( h^{(2)}(X_i, X_j) h^{(2)}(X_k, X_l))| \le CT^{1/3},
\end{equation}
where, $T = [Cov(X_i, X_j)^{\delta} + Cov(X_i, X_k)^{\delta}  + Cov(X_i, X_l)^{\delta}]\times[Cov(X_j, X_i)^{\delta} +\\ Cov(X_j, X_k)^{\delta}  + Cov(X_j, X_l)^{\delta}]\times[Cov(X_k, X_j)^{\delta} + Cov(X_k, X_i)^{\delta}  + Cov(X_k, X_l)^{\delta}]$. Observe that the right-hand side of the inequality in $(3.8)$ has $27$ terms, each term being a product of $3$ covariance terms. Let $r(s) = Cov(X_1, X_{1+s})^ {\gamma}$, $s \ge 0$.  From $(3.8)$, stationarity and observing that $r(0) < \infty$ ($r(j) \le r(0)$, for all $j \ge 1$), 
\begin{align}
  & |E( h^{(2)}(X_i, X_j) h^{(2)}(X_k, X_l))|  \nonumber \\
  & \le C\Big[ r(|i-j|)^2\times[r(|k-i|) + r(|k-j|) + r(|k-l|))]  + r(|i-j|)r(|j-k|)  \nonumber \\
  & \: \: \: + r(|i-j|)r(|j-l|) +  r(|i-k|)r(|j-i|) + r(|i-k|)r(|j-k|) + r(|i-k|)r(|j-l|) \nonumber \\
  &  \: \: \: + r(|i-l|)r(|j-i|) + r(|i-l|)r(|j-k|) + r(|i-l|)r(|j-l|)\Big]. \nonumber \\
    &  \le  C\Big[ r(|i-j|)r(|k-i|) + r(|i-j|)r(|k-j|) + r(|i-j|)r(|k-l|)  \nonumber \\
&    \: \: \: +   r(|i-j|)r(|i-l|) +  r(|i-k|)r(|j-k|) + r(|i-l|)r(|j-k|) +  r(|j-l)r(|j-i|)   \nonumber \\
&   \: \: \: + r(|k-i|)r(|j-l|) + r(|i-l|)r(|j-l|)\Big] = \Delta(i, j, k, l), \text{ (say)}.
 \end{align}
Each term in $(3.9)$ is a product of $2$ distinct covariance terms. It is easy to show that under $\sum_{j=1}^{\infty}r(j) < \infty$, 
\begin{equation}
\sum_{1 \le i < j \le n} \sum_{1 \le k < l \le n} \Delta(i, j, k, l) = o(n^3) \text{ as $n \to \infty$}.
\end{equation}
 Using $(3.10)$, $(3.1)$ follows.
 \end{proof}
\end{lemma}

\begin{remark} The assumption of the function being of bounded Hardy-Krause variation was used in writing the integral in $(3.3)$. Other types of functions could also be considered, provided they generate appropriate measures with bounded total variation.
\end{remark}

\begin{remark} Using Lemmas $2.4$ and $2.5$, Lemma $2.6$ leads to,
\begin{align*}
& |E( h^{(2)}(X_i, X_j) h^{(2)}(X_k, X_l))|  \\
& \le C [Cov(X_i, X_k)^{\delta} + Cov(X_j, X_k)^{\delta} + Cov(X_i, X_l)^{\delta}+ Cov(X_j, X_l)^{\delta}], 
\end{align*}
for some $0 < \delta < 1/2$. Under $\sum_{j=1}^{\infty} Cov(X_1, X_j)^{\delta} < \infty$, 
\begin{equation*}
\sum_{1 \le i < j \le n} \sum_{1 \le k < l \le n} |E( h^{(2)}(X_i, X_j) h^{(2)}(X_k, X_l))| = O(n^3), \text{ as $n \to \infty$}.
\end{equation*}
Hence, Lemma $2.6$ is not enough to ensure $(3.1)$.
\end{remark}

\begin{lemma}
Let the conditions of Lemma $3.1$ be true. Define $U_n$ as the U-statistic based on a symmetric kernel $\rho(x, y)$ which is of  bounded Hardy-Krause variation and left continuous.  Let  $\sigma_1^2$ $=$ $Var(\rho_1(X_1))$  $<$ $\infty$. Define, $\sigma^2_{1j} = Cov( \rho_1(X_1), \rho_1(X_{1+j}))$. Assume $\sum_{j=1}^{\infty}\sigma^2_{1j} < \infty$.
Then,
\begin{equation}
Var(U_n) = \frac{4\sigma^2_U}{n} + o(\frac{1}{n}), \text{ where } \sigma^2_U = \sigma^2_1 + 2\sum_{j=1}^{\infty}\sigma^2_{1j}.
\end{equation}
\begin{proof}
Let $C$ be a generic positive constant in the sequel. Using H-decomposition,
\begin{dmath*}
Var(U_n) = 4Var(H^{(1)}_n) + Var(H^{(2)}_n) + 4Cov(H^{(1)}_n, H^{(2)}_n).
\end{dmath*}
Since $H^{(1)}_n$ $=$ $\frac{1}{n}$ $\sum_{j=1}^{n}h^{(1)}(X_j)$ and $\sum_{j=1}^\infty \sigma_{1j}^2 < \infty$, we get, 
\begin{equation}
Var(H^{(1)}_n) = \frac{1}{n}(\sigma^{2}_{1} +2\sum_{j=1}^{\infty}\sigma^2_{1j})
                          + o(\frac{1}{n}).
\end{equation}
Now, 
\begin{align*}
E(H^{(2)}_{n})^2 & = {{n}\choose{2}}^{-2}\sum_{1 \le i < j\le n}\sum_{1 \le k< l \le n}E\left\lbrace h^{(2)}(X_i, X_j) h^{(2)}(X_k, X_l)\right\rbrace.
 \end{align*}
 As $\rho$ is of bounded Hardy-Krause variation, so is $h^{(2)}$ ($h^{(2)}$ is degenerate by definition). Using Lemma $3.1$,
\begin{equation}
\sum_{1 \le i < j\le n}\sum_{1 \le k< l \le n}E\left\lbrace h^{(2)}(X_i, X_j) h^{(2)}(X_k, X_l)\right\rbrace = o(n^3).
\end{equation}
Therefore,
\begin{equation}
Var(H^{(2)}_{n}) \le E(H^{(2)}_{n})^2  = o(\frac{1}{n}).
\end{equation}
From {$(3.12)$} and {$(3.14)$}, and using Cauchy-Schwarz Inequality we have, 
\begin{align}
|Cov(H^{(1)}_n, H^{(2)}_n)| \le o(n^{-1}).
\end{align}
From {$(3.12)$}, {$(3.14)$} and {$(3.15)$}, we have,
\begin{equation}
Var(U_n) = \frac{4\sigma^2_U}{n} + o(\frac{1}{n}).
\end{equation}
\end{proof}
\end{lemma}

The following gives the central limit theorem for a U-statistic based on a stationary sequence of associated observations with a  kernel of bounded Hardy-Krause variation. 

 \begin{theorem}
Assume the conditions of  Lemma $3.2$ hold and $\sigma^2_U >0$. Suppose there exists a function $\tilde{\rho}_1(\cdot)$ such that ${\rho}_1$ $\ll$ $\tilde{\rho}_1$ and,
\hypertarget{equation 3.11}{
\begin{equation}
\sum_{j=1}^{\infty}Cov(\tilde{\rho}_1(X_1), \tilde{\rho}_1(X_{j})) < \infty.
\end{equation}
}
Then,
\begin{equation}
\frac{\sqrt{n}( U_n - \theta)}{2\sigma_U} \xrightarrow {\mathcal{ L}} N(0, 1) \:\: \text{as}\: \:{n \to \infty}, \text{ where $\sigma^2_U$ is defined by {$(3.11)$}}.
\end{equation}
\begin{proof} Using H-decomposition for $U_n$,
\begin{equation}
\frac{{n^{\frac{1}{2}}}( U_n - \theta)}{2\sigma_U} = n^{{-1}/{2}}\sum_{j=1}^{n}\frac{h^{(1)}(X_j)}{\sigma_U} + n^{1/2}\frac{H^{(2)}_n}{\sigma_U}.
\end{equation}
In addition,
\begin{equation*}
 nE({H^{(2)}_n})^2 \xrightarrow \: 0 \: \text{as}\: \:{n \to \infty}.
\end{equation*}
from $({3.14})$.
Hence,
\begin{equation}
 n^{1/2}\frac{H^{(2)}_n}{\sigma_U}\xrightarrow {p}\: 0 \: \text{as}\: \:{n \to \infty}.
\end{equation}
From {$(3.20)$} and Lemma {2.3}, we get that,
\begin{equation}
n^{{-1}/{2}}\sum_{j=1}^{n}\frac{h^{(1)}(X_j)}{\sigma_U} \xrightarrow {\mathcal{ L}} N(0, 1) \:\: \text{as}\: \:{n \to \infty}.
\end{equation}
Relations {$(3.19)$},  {$(3.20)$}, and  $(3.21)$ prove the theorem.
\end{proof}
\end{theorem} 

\begin{remark} The above results can be easily extended to a U-statistic based on a kernel of any finite degree k. Let $U_n$ be the U-statistic based on a symmetric kernel $\rho(x_1, x_2,\dots, x_k)$ which of finite degree k and of bounded Hardy-Krause variation. Suppose $\sigma_1^2$ $=$ $Var(\rho_1(X_1))$  $<$ $\infty$, $\sum_{j=1}^{\infty}\sigma^2_{1j} <$ $\infty$ and $\sigma^2_U > 0$. Further, let $\sum_{j=1}^\infty Cov(X_1, X_j)^{\gamma}$ $<$ $\infty$, for some $0 < \gamma< 1/6$. Then,
\begin{equation}
Var(U_n) = \frac{k^{2}\sigma^2_U}{n} + o(\frac{1}{n}).
\end{equation}
If the conditions of Theorem $3.3$ hold, then
\begin{equation}
\frac{\sqrt{n}( U_n - \theta)}{k\sigma_U} \xrightarrow {\mathcal{ L}} N(0, 1) \:\: \text{as}\: \:{n \to \infty},
\end{equation}
where $\sigma^2_U$ is defined by {$(3.11)$}.   
\end{remark}

\subsection{Strong law of large numbers}

Christofides $(2004)$ showed that $\lbrace{S_n = {n \choose k}U_n, n \ge k}\rbrace$ ($U_n$ defined by $({1.1})$), is a demimartingale when $E(\rho) = 0$ and $\rho$ is component-wise nondecreasing. Using the concept of demimartingales, he proved a strong law of large numbers for $U_n$ under restrictions on moments of $\rho$. He also extended the result to U-statistics based on kernels $\rho: [a, b] \to \mathbb{R}$ where $[a, b]= [a_1, b_1]\times...\times[a_k, b_k]$ is a $k$-dimensional rectangle and $\rho = h - g$ where, $h$, $g: [a, b] \to \mathbb{R}$ are two component-wise nondecreasing functions and $\Delta_R h \ge 0$ and $\Delta_R g \ge 0$ (given in Definition $1.2$), $ \forall$ $R= [c_1, d_1]\times...\times[c_k, d_k]$ and $a_i \le c_i < d_i \le b_i$ $\forall$ $i =1, 2, ..., k$. 

We observe that kernels which are of bounded Hardy-Krause variation fall into the class of kernels discussed by Christofides $(2004)$. Hence, under restrictions on the moments of the kernel,  as discussed in Theorem $2.1$ and Lemma $2.2$ of Christofides $(2004)$ a strong law of numbers is true for U-statistics based on kernels of bounded Hardy-Krause variation. 

\section{Applications}
\setlength\parindent{24pt}
\subsection{\emph{Gini's  mean difference}}
Suppose we want a measure of variability for observations from a distribution $F$.  A possible index of variability is \emph{Mean difference}, $\theta$, given by,
\begin{equation}
\theta =\int\limits_{{{\mathbb{R}}^2}} |x - y|dF(x)dF(y).
\end{equation}
 Given a sample $\left\lbrace X_j,  1 \le j \le n\right\rbrace$ from $F$, an estimator for $\theta$ is the Gini's mean difference, $U_n$, defined by,
\begin{equation}
U_n = \frac{2}{n(n-1)}\sum_{1 \le i <j \le n} \rho(X_i, X_j), \text{ where the kernel $\rho(x, y) = |x - y|$}.
\end{equation}
 When the observations are $i.i.d$, using Hoeffding $(1948)$ it can be shown that $\frac{\sqrt{n}(U_n - \theta)}{2\sqrt{\mathfrak{F}- \theta^2}}  \xrightarrow {\mathcal{ L}} N(0, 1)$ as ${n \to \infty}$, where $\mathfrak{F} = \int\limits_{{{\mathbb{R}}^3}} |x-y||x-z|dF(x)dF(y)dF(z)$.

We now obtain the limiting distribution of $U_n$ when the observations are stationary and associated.
\begin{theorem} Let $\{X_n, n\ge 1\}$ be a sequence of stationary associated random variables having one dimensional marginal distribution function $F$ and a bounded density function. Assume $|X_n| < C_1 < \infty$, for all $n \ge 1$. Let,
\begin{equation}
 \sum_{j=1}^{\infty}Cov(X_1, X_j)^{\gamma} \: < \: \infty, \text{ for some $0 < \gamma < 1/6$}.
\end{equation}
Then,
\begin{equation*}
\frac{\sqrt{n}(U_n - \theta)}{2\sigma_U}  \xrightarrow {\mathcal{ L}} N(0, 1) \:\: \text{as}\:\: {n \to \infty},
\end{equation*}
 where, $U_n$ and $\sigma^2_U$  are defined by  $({4.2})$ and $({3.11})$ respectively.
\begin{proof}
By H-decomposition, $U_n = \theta + 2H^{(1)}_n + H^{(2)}_n$, where $\theta$ is given by $({4.1})$.
We observe $\rho(x, y)$ is a function of bounded Hardy-Krause variation (Definition ${1.3}$).
Now, $\rho_1(x) = \int_{-\infty}^{\infty}|x - y|dF(y)$. We can choose $\tilde{\rho}_1(x) = Cx$, for some $C > 0$ as $\rho_1(\cdot)$ is Lipschitzian. 
Using $({4.3})$ and Theorem $3.3$, we have
\begin{equation*}
\frac{\sqrt{n}(U_n - \theta)}{2\sigma_U} \xrightarrow {\mathcal{ L}} N(0, 1) \:\: \text{as}\: \:{n \to \infty}.
\end{equation*}
Here, $  Var(\rho_1(X_1)) = \mathfrak{F} - \theta^2$ and $Cov(\rho_1(X_1), \rho_1(X_j)) = E[|X_1 - X||X_j - Y|] - \theta^2$, where $X$ and $Y$ are independent of $\lbrace{X_n, n \ge 1}\rbrace$ and are  $i.i.d$ such that $F$ is the marginal distribution function of $X$.
\end{proof}
\end{theorem}

\begin{remark}The above result can be extended to random variables which are not uniformly bounded using the usual truncation techniques.
\end{remark}

The next result is needed in simulation analysis. For details, see comment $(4)$ in section 5.
 
\begin{theorem}  Let $\{X_n, n\ge 1\}$ be a sequence of stationary associated random variables having one dimensional marginal distribution function $F$. Let $|X_n| < C_1 < \infty$, for all $n \ge 1$. Assume,
 \hypertarget{equation 4.4}{
 \begin{equation}
 \sum_{j=n+1}^{\infty}Cov(X_1, X_j) = O(n^{-(k-2)/2}), \text{ for some $k>2$, $k \in \mathbb{N}$. }
\end{equation}}
\begin{equation*}
\text{Then, }\underset{x \in  [-C_1, C_1]\: }{\: sup}\Big|  \frac{\sum_{j=1}^n (|X_j - x| - E|X_j -x|)}{b_n}\Big| \to 0 \text{ a.s as n } \to \infty, \text{ where,}
\end{equation*}
$b_n = O(n^{1+u/2-p})$, for some $u > 1$ and $p \in (0, 1)$ such that $\frac{k}{2}(1+u) - (k+1)p > 1$. 
\begin{proof} Divide the interval $ [-C_1, C_1]$ into $d_n = O(n^p)$, for $p \in (0, 1)$ small intervals as follows:
Let $-C_1 = y_{n_0} <  y_{n_1} <  ... <  y_{n_{d_n}} = C_1$. The $d_n$ intervals are denoted as $I_{n_i} = [ y_{n_{i-1}}, y_{n_i}], i = 1, 2, ..., d_n$, each of length $\delta_n = \frac{2C_1}{d_n}$. Let $x_{n_i} \in I_{n_i}$.
 \begin{align}
& \underset{x \in  [-C_1, C_1]\: }{sup}\Big|  \frac{\sum_{j=1}^n (|X_j - x| - E(|X_j -x|) )}{b_n}\Big| = \underset{i}{max}\underset{x \in I_{n_i}}{sup} \Big|  \frac{\sum_{j=1}^n (|X_j - x| - E(|X_j -x|) )}{b_n}\Big| \nonumber \\
&  \le \underset{i}{max}\underset{x \in I_{n_i}}{sup} \Big|  \frac{\sum_{j=1}^n (|X_j - x| - |X_j -x_{n_i}| )}{b_n}\Big| +  \underset{i}{max}\underset{x \in I_{n_i}}{sup}\Big|  \frac{\sum_{j=1}^n (|X_j - x_{n_i}| - E(|X_j -x_{n_i}|) )}{b_n}\Big|  \nonumber \\
& + \underset{i}{max}\underset{x \in I_{n_i}}{sup} \Big|  \frac{\sum_{j=1}^n (E(|X_j - x|) - E(|X_j -x_{n_i}|) )}{b_n}\Big|  = I_1 + I_2 + I_3 \text{ (say)}. 
\end{align}
For ${x \in I_{n_i}}$, $\Big|  \frac{\sum_{j=1}^n (|X_j - x| - |X_j -x_{n_i}| )}{b_n}\Big|$ $\le$ $\frac{ n\delta_n}{b_n}$. Hence,  $I_1 \le  \frac{n\delta_n}{b_n}$. Similarly, $I_3 \le \frac{n\delta_n}{b_n}$.\\
$I_2 = \underset{i}{max}\underset{x \in I_{n_i}}{sup}\Big|  \frac{\sum_{j=1}^n (|X_j - x_{n_i}| - E(|X_j -x_{n_i}|) )}{b_n}\Big| = \underset{i}{max}\Big|  \frac{\sum_{j=1}^n (|X_j - x_{n_i}| - E(|X_j -x_{n_i}|) )}{b_n}\Big| $. \\
For any $\epsilon > 0$, 
 \begin{align}
& P\Big[\underset{i}{max}\Big|  \frac{\sum_{j=1}^n (|X_j - x_{n_i}| - E(|X_j -x_{n_i}|) )}{b_n}\Big|  > \epsilon \Big] \nonumber \\
& \le \sum_i  P\Big[ \Big|  \frac{\sum_{j=1}^n (|X_j - x_{n_i}| - E(|X_j -x_{n_i}|) )}{b_n}\Big|  > \epsilon \Big] \nonumber \\
& = \sum_i  P\Big[ \Big|  \frac{S_n(x_{n_i} )}{b_n}\Big|  > \epsilon \Big] \le d_n \underset{i}{max}\frac{E |S_n(x_{n_i} )|^k}{b_n^k\epsilon^k}. 
\end{align}
where $S_n(x_{n_i}) = \sum_{j=1}^n (|X_j - x_{n_i}| - E|X_j -x_{n_i}|)$. Let $B$ be a generic positive constant in the sequel.
 \begin{equation}
 \frac{E |S_n(x_{n_i} )|^k}{b_n^k \epsilon^k }  =  \frac{1}{2^k}\frac{E |S_n(x_{n_i}) - \tilde{S}_n(x_{n_i})+ S_n(x_{n_i}) + \tilde{S}_n(x_{n_i})|^k}{b_n^k \epsilon^k},
 \end{equation}
   where, $\tilde{S}_n(x_{n_i}) = \sum_{j=1}^n B(X_j - x_{n_i})$. Observe, $|X_ j - x| \ll  B(X_j -x)$, $j \ge 1$, for all $x \in [-C_1, C_1]$. For all $x \in [-C_1, C_1]$,
  \begin{align*}
  & \sum_{j= n+1}^{\infty} Cov(|X_1 - x| - E|X_1 -x| - B(X_1 -x) , |X_j - x| - E|X_j -x| - B(X_j -x))  \nonumber \\
& \le B\sum_{j=n+1}^{\infty} Cov(X_1 -x, X_j -x) =  B\sum_{j=n+1}^{\infty}Cov(X_1, X_j) = O(n^{-(k-2)/2}).
  \end{align*}
  Similarly, $\sum_{j= n+1}^{\infty} Cov(|X_1 - x| - E|X_1 -x| + B(X_1 -x), |X_j - x| - E|X_j -x|+ B(X_j -x)) = O(n^{-(k-2)/2})$, for all $ x \in [-C_1, C_1]$. 
  
  Using the equality $(a+b)^k = \sum_{j=0}^k {k \choose j} a^j b^{k-j}$ and the Cauchy-Shwarz inequality, we have, $E |S_n(x_{n_i}) - \tilde{S}_n(x_{n_i})+ S_n(x_{n_i}) + \tilde{S}_n(x_{n_i})|^k$  $\le \sum_{j=0}^k { k \choose j} (E|S_n(x_{n_i}) - (\tilde{S}_n(x_{n_i})- E(\tilde{S}_n(x_{n_i})))|^{2j}E|S_n(x_{n_i}) + (\tilde{S}_n(x_{n_i})- E(\tilde{S}_n(x_{n_i})))|^{2(k-j)})^{1/2}$.  
  
 Observe that both $\lbrace{ |X_j- x_{n_i}| - E|X_j - x_{n_i}| - B(X_j -x_{n_i}) ; j \ge 1}\rbrace$  and \\ $\lbrace{ |X_j - x_{n_i}| - E|X_j - x_{n_i}| + B(X_j -x_{n_i}) ; j \ge 1}\rbrace$ form an associated sequence. Hence, using Lemma $2.8$,  we get, $ \frac{E |S_n(x_{n_i} )|^k}{b_n^k\epsilon^k} \le B \frac{(n^{j} n^{k-j})^{1/2}}{b_n^k\epsilon^k 2^k}$.
 
 From $({4.5})$, we have,
\begin{align}
& \underset{x \in  [-C_1, C_1]}{sup} \Big| \frac{\sum_{j=1}^n (|X_j -x| - E|X_j -x|)}{b_n} \Big| \le \frac{2n\delta_n}{b_n} + \underset{i}{max}\Big| \frac{\sum_{j=1}^n (|X_j -x_{n_i}| - E|X_j -x_{n_i}|)}{b_n} \Big|.
\end{align}
 \begin{align}
& P\Big[  \underset{x \in  [-C_1, C_1]}{sup} \Big| \frac{\sum_{j=1}^n (|X_j -x| - E|X_j -x|)}{b_n} \Big| > \epsilon \Big] \nonumber \\
& \le P\Big[ \frac{2n\delta_n}{b_n} > \frac{\epsilon}{2} \Big] + P \Big[ \underset{i}{max}\Big| \frac{\sum_{j=1}^n (|X_j -x_{n_i}| - E|X_j -x_{n_i}|)}{b_n} \Big| > \frac{\epsilon}{2} \Big]  \nonumber \\
& \le P\Big[ \frac{4C_1n}{b_n d_n} > \frac{\epsilon}{2} \Big] + B \frac{n^{k/2}d_n}{b_n^k\epsilon^k} \le \Big(\frac{8C_1n}{b_nd_n\epsilon}\Big)^2 + B \frac{n^{k/2}d_n}{b_n^k\epsilon^k}.
\end{align}
Finally,
  \hypertarget{equation 4.12}
  {\begin{align}
& \sum_{n = 1}^{\infty}P\Big[  \underset{x \in  [-C_1, C_1]}{sup} \Big| \frac{\sum_{j=1}^n (|X_j -x| - E|X_j -x|)}{b_n} \Big| > \epsilon \Big] 
 \le B\sum_{n = 1}^{\infty}  \Big\lbrace \Big(\frac{n}{b_nd_n}\Big)^2 + \frac{n^{k/2}d_n}{b_n^k \epsilon^k} \Big\rbrace.
\end{align}}
Result follows if $b_n = O(n^{1+u/2-p})$, for some $u > 1$ and $\frac{k}{2}(1+u) - (k+1)p > 1$. 
\end{proof}
\end{theorem}

\subsection{Empirical joint distribution functions} 

Let $\lbrace{X_n , n \ge 1}\rbrace$ be a sequence of stationary associated random variables. Bagai and Prakasa Rao $(1991)$ had discussed the asymptotics for the empirical estimator of  survival function for this sequence. Henriques and Oliveira $(2003)$  had discussed the asymptotics for the histogram estimator for the two-dimensional distribution function of $(X_1, X_{k+1})$. In both the cases, the kernel is of bounded Hardy-Krause variation. Similarly, the kernel of the histogram estimator for any finite k-dimensional distribution function is also of bounded Hardy-Krause variation. Assume further $|X_n| < C_1 , n \ge 1$, $0 < C_1 < \infty$ . Then, for the two-dimensional distribution function of $(X_1, X_{k+1})$, $P(X_1 \le s, X_{k+1} \le t)$, $s, t \in \mathbb{R}$, the histogram estimator is,
\begin{equation}
U_n(s, t)= \frac{1}{n-k}\sum_{i=1}^{n-k} Y_{i, i+k}
\end{equation}
 where, $Y_{i, i+k} = I(X_i \le s, X_{i+k} \le t)$. Using Lemmas $2.4$, $2.5$ and $2.6$,
 \begin{align*}
 Cov(Y_{i, i+k}, Y_{j, j+k})  \le & \: C(Cov(X_i, X_j)^{\delta} + Cov(X_{i+k}, X_j)^{\delta} + Cov(X_i, X_{j+k})^{\delta} + \nonumber \\
&  Cov(X_{i+k}, X_{j+k})^{\delta}), 
 \end{align*}
 for some $0 < \delta < 1/2$. Using Lemma $2.1$ and the condition $\sum_{j=1}^{\infty} Cov(X_1, X_j)^{\delta} < \infty$ the asymptotic normality $U_n(s, t)$ can be obtained.
  \begin{remark} Henriques and Oliveira $(2003)$ had obtained the asymptotic distribution of the histogram estimator for stationary associated random variables belonging to $L_2$ under the covariance restriction of $\sum_{j=1}^{\infty} Cov(X_1, X_j)^{1/3} < \infty$. Replacing the covariance control by Demichev's inequality (Lemma $2.5$) leads to the condition on the covariances being less restrictive. 
 \end{remark}

\section{Simulation Analysis}

 The asymptotic normality of Gini's mean difference based on stationary and associated observations are investigated via simulations. Let  $\lbrace {Y_j, j \ge 1} \rbrace$ be $i.i.d$  from $Exp(1/m)$ for some $m \in \mathbb{N}$. If $X_j$ $=$ $min(Y_j,\cdots,Y_{j+m-1})$ for all $1 \le j \le n$, then $\lbrace{X_j, 1 \le j \le n}\rbrace$ forms a set of stationary associated random variables such that $X_j's$ are standard exponential variables $(Exp(1))$. Similarly, in order to obtain stationary associated random variables $\lbrace{X_j,  1 \le j \le n}\rbrace$ such that $X_j's$ are standard normal variables $(N(0, 1))$, we can set $X_j = Y_j+  \cdots+ Y_{j+m-1}$, where $\lbrace {Y_j,  j \ge 1} \rbrace$ are $i.i.d$ from $N(0,1/m)$ for some $m \in \mathbb{N}$.
\begin{itemize}
\item[(1)] We use the statistical software R (http://www.r-project.org; R Development Core Team $(2011)$) for our simulations. The samples $\lbrace{X_j, 1 \le j \le n}\rbrace$ are generated as follows.
\begin{itemize}
\item [(S1)]  $X_j = min(Y_j, Y_{j+1})$, where $\lbrace {Y_j, j \ge 1} \rbrace$ are pseudo-random numbers from $Exp(1/2)$ generated using $rexp$ function in R.
 \item[(S2)]  $X_j =min( Y_j, Y_{j+1}, Y_{j+2})$, where $\lbrace {Y_j, j \ge 1} \rbrace$ are pseudo-random numbers from $Exp(1/3)$ generated using $rexp$ function in R.
 \item [(S3)]  $X_j = min(Y_j, \cdots, Y_{j+9})$, where $\lbrace {Y_j, j \ge 1} \rbrace$ are pseudo-random numbers from  $Exp(1/10)$ generated using $rexp$ function in R.
\item [(S4)] $X_j= Y_j  + Y_{j+1}$, where $\lbrace {Y_j, j \ge 1} \rbrace$ are pseudo-random numbers from $N(0, 1/2)$, generated using $rnorm$ function in R. 
\item[(S5)] $X_j=  Y_j +  Y_{j+1}+ Y_{j+2}$, where $\lbrace {Y_j, j \ge 1} \rbrace$ are pseudo-random numbers from  $N(0, 1/3)$, generated using $rnorm$ function in R. 
\item[(S6)] $X_j=  Y_j + \cdots + Y_{j+9}$, where $\lbrace {Y_j, j \ge 1} \rbrace$ are pseudo-random numbers from $N(0, 1/10)$, generated using $rnorm$ function in R. 
\end{itemize}
\item [(2)] The results are based on 10,000 replications and $\alpha$ = 0.05.
\item [(3)] For our simulations, we  use Lemma $2.7$ for the estimation of $\sigma_U$. We choose $\ell_n = [n^{3/5}]$, smallest integer less than or equal to $n^{3/5}$.
\hypertarget{pt5}{\item [(4)]  In Lemma $2.7$, $Y_i = \rho_1(X_i)$ and $\tilde{Y}_i = bX_i$, for some constant $b$, $b > 0$ $\forall$ $i \ge 1$, as $\rho_1$ is lipshitz.  For practical applications, the distribution function of the underlying population $F$ will be unknown. Hence, an estimator for $\rho_1(x)$ is needed. Let $\hat{B}_n$  be analogous to $B_n$ with $S_j (k)$ replaced by $\hat{S}_j(k) = \sum_{i = j+1}^{j+ k}\hat{\rho}_1(X_i)$, and $\bar{Y}_n$ by $\bar{\hat{Y}}_n = \sum_{i =1}^n{\hat{\rho}_1(X_i)}$, where $\hat{\rho}_1(x) = \sum_{j=1}^n \frac{|X_j - x|}{n}$. Define, $|Z_i| =$ $2$ $|\hat{\rho}_1(X_i) - \rho_1(X_i)|$.
\begin{align*}
|B_n - \hat{B}_n| & = |\frac{1}{n-\ell+1} \sum_{j=1}^{n-\ell} \frac{|S_j(\ell) - \ell \bar{Y}_n|}{\sqrt{\ell}} - \frac{1}{n-\ell+1} \sum_{j=1}^{n-\ell} \frac{|\hat{S}_j(\ell) - \ell \bar{\hat{Y}}_n|}{\sqrt{\ell}}| \nonumber \\
& \le \frac{1}{(n-\ell+1)\sqrt{\ell}}  \sum_{j=1}^{n-\ell} \Big( \sum_{i =j +1}^{j +\ell}|Z_i| + \frac{\ell}{n} \sum_{i=1}^n|Z_i| \Big) \nonumber \\
& \le \frac{2}{(n-\ell+1)\sqrt{\ell}} (n - \ell) \ell  \: \underset{x}{ sup} \: |\hat{\rho}_1(x) - \rho_1(x)| \\
  & = 2  \: \frac{\sqrt{\ell}}{n^s} \: \underset{x}{sup}\: | \sum_{j=1}^n \frac{|X_j - x| - \rho_1(x)}{n^{1-s}}|
\end{align*}
Putting $s = \frac{3}{10}$, we get $\frac{\sqrt{\ell}}{n^s} = O(1)$ ($\ell = \ell_n = [n^{3/5}]$).  In Theorem $4.2$, putting $k=10$, $p = 17/20$, $u = 11/10$, and assuming $\sum_{j=n+1}^{\infty} Cov(X_1, X_j) = O(n^{-4})$, we get, $\underset{x}{sup}| \sum_{j=1}^n \frac{|X_j - x| - \rho_1(x)}{n^{7/10}}| \to 0$ a.s as $n \to \infty$. Hence, $|B_n - \hat{B}_n|$ $\to 0$ a.s as $n \to \infty$.}
\end{itemize}

In the following tables, 
\begin{itemize}
\item [a)]$\bar{g}$ denotes the mean of the $r = 10,000$ sample Gini's  mean difference values, $g_i$, $1 \le i \le r$; 
\item [b)]  E.M.S.E (g) $= \frac{1}{r -1} \sum_{i=1}^{r}(g_i - \bar{g})^2$, where $E.M.S.E$ denotes Estimated M.S.E;
\item [c)] C.P (g) $= \frac{N}{r}$, where $N$ $=$ $\# $ $\lbrace{i : g_i \in (\bar{g} - 2\bar{\hat{B}}_n\times\frac{z_{0.025}}{\sqrt{n}}, \bar{g} + 2\bar{\hat{B}}_n\times\frac{z_{0.025}}{\sqrt{n}} )}\rbrace$. Here, $\bar{\hat{B}}_n = \frac{1}{r}\sum_{i=1}^r \hat{B}_n(i)$, $z_{0.025} = 1.959964$, and ${\hat{B}}_n(i)$, $1 \le i \le r$, denotes the estimated value for each sample;
\item [d)] E.M.S.E (${\bar{{B}}_n}$) $= \frac{1}{r -1} \sum_{i=1}^{r}(\hat{B}_n(i) - \bar{\hat{B}}_n)^2$; 
\item [e)] Median (g), Skewness (g), and Kurtosis (g)  are the corresponding characteristics of the $r$ sample statistic values.
\end{itemize}

\begin{center}
{\scriptsize\begin{tabular}{p{0.05in}p{.6in}  p{.6in}p{.6in} p{0.6in}  p{.6in} p{.6in}| }
 \multicolumn{7}{c}{\textbf{ Table 5.1}\: \emph{\textbf{Simulation Results for $Exp (1)$}}}\\
 \cline{1-7}
  \multicolumn{1}{|c|}{ $(S1)$ (m=2)} & n=50 & n=100 & n=200  & n=300 & n=500 & n=1000\\ \hline
\multicolumn{1}{|c|}{$\bar{g}$}&  0.9836181 & 0.9930083 & 0.9980009 & 0.9983753 &  0.9984739 &  0.9989455\\
\multicolumn{1}{|c|}{E.M.S.E (g)} &  0.03676916 & 0.01931296 & 0.009705885 & 0.006434793 & 0.003937328 & 0.001930927\\
\multicolumn{1}{|c|}{C.P (g)} & 0.8974 & 0.9133 & 0.9258 & 0.9316 & 0.9342 & 0.9403\\
\multicolumn{1}{|c|}{Median (g) } &  0.9711876 & 0.9847915 & 0.9941595 &   0.9947342 & 0.9978153 & 0.9982964\\ 
\multicolumn{1}{|c|}{Skewness (g)} & 0.3965101 & 0.3466176 &  0.2007335 &  0.200636 & 0.1203145 & 0.08936168\\  
\multicolumn{1}{|c|}{Kurtosis (g) } &   3.256429 & 3.219073 & 3.026955 & 3.002957 & 3.001628 & 3.020152\\ \hline  \\ \hline
  \multicolumn{1}{|c|}{  $(S2)$ (m=3)} & n=50 & n=100 & n=200  & n=300 & n=500 & n=1000\\ \hline
\multicolumn{1}{|c|}{$\bar{g}$} & 0.9732993 & 0.9875498 & 0.9941661 &  0.9972799 & 0.9972273 & 0.9980701\\
\multicolumn{1}{|c|}{E.M.S.E (g)} &  0.04972117 & 0.02671703 & 0.0133838 & 0.008855159 & 0.00531166 & 0.002623486 \\
\multicolumn{1}{|c|}{C.P (g)} & 0.8804 & 0.9002 & 0.9201 & 0.9241 & 0.9296 & 0.9396\\
\multicolumn{1}{|c|}{Median (g)} & 0.9563219 & 0.9793019 & 0.9885436 & 0.9928442 & 0.9950708  &  0.9966722\\ 
\multicolumn{1}{|c|}{Skewness (g) }& 0.5218394 & 0.3517759 & 0.2756275 & 0.2407178  &0.1604277  &  0.1284387\\  
\multicolumn{1}{|c|}{Kurtosis (g) }& 3.513031 & 3.167723 & 3.106484 & 3.15176 & 3.057265 & 3.060678\\ \hline \\ \hline
 \multicolumn{1}{|c|}{ $(S3)$ (m=10)}& n=50 & n=100 & n=200  & n=300 & n=500 & n=1000\\ \hline
\multicolumn{1}{|c|}{$\bar{g}$}& 0.8758573 & 0.93135 & 0.9627473 & 0.9779133 &  0.9869055 & 0.9940865\\
\multicolumn{1}{|c|}{E.M.S.E (g)} &  0.1280773 & 0.07158817 & 0.03887426 &  0.02633781 & 0.0164874 &  0.008295358 \\
\multicolumn{1}{|c|}{C.P (g)}& 0.7579 &   0.8134 & 0.8478 & 0.8723 & 0.8889 & 0.9097 \\
\multicolumn{1}{|c|}{Median (g) } & 0.8260538 & 0.9033452 & 0.949436 &  0.9671378 & 0.9795772 & 0.9908148\\ 
\multicolumn{1}{|c|}{Skewness (g)} & 0.9151533 &  0.6293821 & 0.4272413 & 0.3993889 & 0.3150814 &  0.2282979\\  
\multicolumn{1}{|c|}{Kurtosis (g) } &  4.462582 & 3.537275 & 3.239824 & 3.336504 & 3.231984 & 3.017788\\ \hline 
 \end{tabular}}
\end{center}

\begin{center}
{\scriptsize\begin{tabular}{p{0.05in}p{.6in}  p{.6in}p{.6in} p{0.6in}  p{.6in} p{.6in}| }
 \multicolumn{7}{c}{\textbf{ Table 5.2}\: \emph{\textbf{Simulation Results for $N (0, 1)$}}}\\
 \cline{1-7}
  \multicolumn{1}{|c|}{$(S4)$ (m = 2)}& n=50 & n=100 & n=200  & n=300 & n=500 & n=1000\\ \hline
\multicolumn{1}{|c|}{$\bar{g}$}&   1.114259 & 1.120531 & 1.125395  & 1.126033 & 1.127356 & 1.127414\\
\multicolumn{1}{|c|}{E.M.S.E (g)} &  0.01965189 & 0.00962333 & 0.004895764 & 0.003354603 & 0.001963522 & 0.0009780631\\
\multicolumn{1}{|c|}{C.P (g)}& 0.9077 &  0.9262 & 0.9317 & 0.9327& 0.9413 & 0.9441\\
\multicolumn{1}{|c|}{Median (g)}& 1.111408 & 1.119302 & 1.124416 &  1.124914 & 1.127125 & 1.127346\\
\multicolumn{1}{|c|}{Skewness (g) } &  0.1287656 & 0.1264063 & 0.09834354 & 0.06234944 &   0.07745284 & 0.04384436 \\  
\multicolumn{1}{|c|}{Kurtosis (g) } &  2.993197 & 2.923884 & 2.967621 & 2.907618 & 2.947061 & 2.880582 \\ \hline  \\ \hline
  \multicolumn{1}{|c|}{$(S5)$ (m =3)} & n=50 & n=100 & n=200  & n=300 & n=500 & n=1000\\ \hline
\multicolumn{1}{|c|}{$\bar{g}$} &  1.102877 & 1.114272 &  1.122001 & 1.122333 & 1.125143 & 1.126909\\
\multicolumn{1}{|c|}{E.M.S.E (g)}&  0.02673821 & 0.01331181 & 0.006773474 &  0.004577988  &  0.002688516 & 0.001317823\\
\multicolumn{1}{|c|}{C.P (g)} & 0.8923 & 0.9173 & 0.9259 & 0.9317 & 0.9425 & 0.9483 \\
\multicolumn{1}{|c|}{Median (g) }&  1.096844 & 1.11083 & 1.119453 & 1.121584 &  1.124423 &  1.126417\\ 
\multicolumn{1}{|c|}{Skewness (g) } & 0.2345972 &  0.1703191 & 0.1451524 &  0.08957866 & 0.05842091&  0.04335457  \\  
\multicolumn{1}{|c|}{Kurtosis (g) } & 3.027241 &  3.073061 & 2.992162 & 3.030902 & 2.968134 & 2.987061 \\ \hline \\ \hline
 \multicolumn{1}{|c|}{$(S6)$ (m=10)} & n=50 & n=100 & n=200  & n=300 & n=500 & n=1000\\ \hline
\multicolumn{1}{|c|}{$\bar{g}$}& 1.004597 & 1.064742 & 1.096209 & 1.10654 & 1.113758 & 1.121415 \\
\multicolumn{1}{|c|}{E.M.S.E (g)}& 0.06772895 & 0.03811105 & 0.02010571 &  0.01393415 & 0.008248388 & 0.004292068  \\
\multicolumn{1}{|c|}{C.P (g)}&  0.7678 & 0.8503 & 0.8878 & 0.8992 & 0.9186 & 0.9277\\
\multicolumn{1}{|c|}{Median(g) }&  0.9787768 & 1.051586 & 1.088681 & 1.101294 &  1.112035 &  1.119691\\ 
\multicolumn{1}{|c|}{Skewness(g) }& 0.6010408 & 0.4048426 & 0.2736039 & 0.2290044 & 0.2082488  &  0.1407511  \\  
\multicolumn{1}{|c|}{Kurtosis(g)}& 3.470072 &  3.252502 & 3.034935 & 3.044965 & 3.075823 &  3.000726\\ \hline 
 \end{tabular}}
\end{center}

\underline{{\textbf{Observations}}}
\begin{itemize}
\item[(i)] \emph{Estimation of $\sigma_U$}: As discussed earlier, we use an estimator for $\sigma_U$ for simulations. {$(4)$} and Lemma $2.7$  imply that $\sqrt{{\pi}/{2}}\hat{B}_n$ is also a consistent estimator for $\sigma_U$. For the sample generated from $Exp (1)$, using $(S1)$, $(S2)$,  and $(S3)$, we analyze the performance of the estimator by comparing $2\sqrt{{\pi}/{2}}\hat{B}_n$ with the actual values ($2\sigma_U$).  The following table shows that as the sample size increases, the value of bias reduces. As expected, $E.M.S.E$ (Estimated M.S.E) also reduces with the increase in the sample size. For $m =2, 3$, the rate of convergence is faster than for $m = 10$.
\begin{center}
{\scriptsize\begin{tabular}{p{0.05in}|p{.4in}  p{.4in}p{.4in} p{0.4in}  p{.4in} p{.4in}| }
 \multicolumn{7}{c}{\textbf{ Table 5.3}\: \emph{\textbf{{Performance of $\hat{B}_n$ for Exp(1)}}} }\\
 \cline{1-7}
  \multicolumn{1}{|c|}{ $(S1)$ (m=2), $2\sigma_U = 1.393864$ } & n=50 & n=100 & n=200  & n=300 & n=500 & n=1000\\ 
\multicolumn{1}{|c|}{ $2\sqrt{{\pi}/{2}}\bar{\hat{B}}_n$ }& 1.113067 & 1.206764 &  1.26211 & 1.285339 & 1.312719  & 1.335088 \\
\multicolumn{1}{|c|}{ Bias = 2$|\sqrt{{\pi}/{2}}\bar{\hat{B}}_n - \sigma_U|$ }& 0.280797 & 0.187100 & 0.131754 & 0.108525 & 0.081145 & 0.058776 \\
\multicolumn{1}{|c|}{E.M.S.E (2$\sqrt{{\pi}/{2}}{\hat{B}}_n$)}&  0.1941304 & 0.1271084 & 0.0828345 & 0.0641774 &  0.0455669 & 0.029956\\ \hline
 \multicolumn{1}{|c|}{  $(S2)$ (m=3),  $2\sigma_U = 1.639871$} & n=50 & n=100 & n=200  & n=300 & n=500 & n=1000\\ 
\multicolumn{1}{|c|}{$2\sqrt{{\pi}/{2}}\bar{\hat{B}}_n$} & 1.217808 & 1.345743 & 1.434806 & 1.470134 & 1.506485  & 1.544193\\
\multicolumn{1}{|c|}{ Bias = 2$|\sqrt{{\pi}/{2}}\bar{\hat{B}}_n - \sigma_U|$ }& 0.422063 & 0.294128 & 0.205065 & 0.169737  & 0.133386 & 0.095678 \\
\multicolumn{1}{|c|}{E.M.S.E ($2\sqrt{{\pi}/{2}}{\hat{B}}_n$)}&  0.2784683 & 0.1971176 & 0.1303082 & 0.1002306 &1.508937 &  0.046944\\ \hline
 \multicolumn{1}{|c|}{  $(S3)$ (m=10), $2\sigma_U = 2.897561$} & n=50 & n=100 & n=200  & n=300 & n=500 & n=1000\\ 
\multicolumn{1}{|c|}{$2\sqrt{{\pi}/{2}}\bar{\hat{B}}_n$} & 1.435125 & 1.761425 & 2.011354 & 2.149753 & 2.29772 & 2.469504 \\
\multicolumn{1}{|c|}{ Bias = 2$|\sqrt{{\pi}/{2}}\bar{\hat{B}}_n - \sigma_U|$  }&  1.462436 &  1.136136 & 0.886207 & 0.747808 & 0.599841 & 0.428057 \\
\multicolumn{1}{|c|}{E.M.S.E ($2\sqrt{{\pi}/{2}}{\hat{B}}_n$)}&  0.7394865 & 0.6412837 & 0.4948429 & 0.4134175 & 0.3136523 & 0.205857  \\ \hline
 \end{tabular}}
\end{center}
\item[(ii)] \emph{Asymptotic Normality}: From Tables $5.1$ and $5.2$, we observe that for a fixed $m$ as the sample size increases, the approximation to the normal distribution is better.  For $m =2, 3$, the convergence to normality is faster, as expected, as the variables are ``almost independent''. For $m =10$, we see that the approximation is good only for much  larger values of $n$. The use of the estimator  of $\sigma_U$ could also affect the convergence as the bias and $E.M.S.E$ (Estimated M.S.E)  reduces much faster for $m = 2, 3$ than for $m = 10$. 
\item[(iii)] \emph{Estimation of the mean difference}: When $X_j's$ are $Exp(1)$, the value of the mean difference, $\theta$, is $1$. From Table {5.1}, it can be seen that when $m = 2, 3$, the convergence of the mean of 10,000 sample Gini's mean difference values to 1 is faster than when $m = 10$. This is expected as greater dependence leads to a slower rate. Similar results are observed from Table {5.2}. Here $\theta = 1.128379$.
\item[(iv)] \emph{Comparison with $i.i.d$ setup}: A comparison of the simulation results with the results of Greselin and Zenga $(2006)$ who had performed the simulations for the statistic under the $i.i.d$ setup, indicate that relatively larger sample sizes are needed for applying the asymptotic results under the dependent setup than under the $i.i.d$ setup. 
\end{itemize}

\section{Conclusions}
In this paper, we give the limiting distribution of U-statistics based on kernels of bounded Hardy-Krause variation when the underlying sample consists of stationary associated observations. As an application, we obtain the asymptotic distribution of Gini's mean difference under the dependent setup. Simulation results performed for the statistic indicate that reasonable sample sizes are needed for using the normality approximation. Greater the dependence, larger the sample sizes needed for a viable use of the asymptotic normality results. 

 Results for kernels which are differentiable have been discussed in Garg and Dewan $(2015)$. Results for discontinuous kernels that are not component-wise monotonic and are not functions of bounded Hardy-Krause variations are under preparation.
 
 \section{Acknowledgment }
We would like to thank the two anonymous referees for their critical comments that helped in improving the earlier version of the paper.

\bibliographystyle{yearfirst}

{\footnotesize{

\end{document}